\documentclass[thmsa]{article}
\usepackage{amssymb}
\usepackage{amsfonts}
\usepackage{graphicx}

\newtheorem{theorem}{Theorem}

\newtheorem{lemma}[theorem]{Lemma}

\newtheorem{remark}[theorem]{Remark}

\def\proof{\noindent{{\bf Proof. }}}
\def\endproof{ $\blacksquare$}
\def\dfrac{\displaystyle\frac}

\begin{document}

\title{A nonexistence result for a nonlinear elliptic equation with singular and decaying potential}

\author{{\small \smallskip }Marino Badiale\thanks{%
Partially supported by the PRIN2009 grant ``Critical Point Theory and
Perturbative Methods for Nonlinear Differential Equations''}~ - Michela
Guida - Sergio Rolando\footnotemark[1]  \smallskip\\
\textit{\small Dipartimento di Matematica}\\
\textit{{\small Universit\`{a} degli Studi di Torino}}, \textit{\small Via
Carlo Alberto 10, 10123 Torino}, \textit{\small Italy}\\
\textit{\small e-mail: }{\small marino.badiale@unito.it}\textit{\small , }%
{\small michela.guida@unito.it}\textit{\small , }{\small %
sergio.rolando@unito.it}}
\date{}
\maketitle

\begin{abstract}
Several existence and nonexistence results are known for positive solutions 
$u\in D^{1,2}(\mathbb{R}^{N})\newline
\cap L^{2}(\mathbb{R}^{N},\left| x\right| ^{-\alpha}dx)\cap L^{p}(\mathbb{R}^{N})$
to the equation 
\[
-\triangle u+\frac{A}{\left| x\right| ^{\alpha }}u=u^{p-1}\quad %
\textrm{in }\mathbb{R}^{N}\setminus \left\{ 0\right\} ,\quad N\geq 3,~A,\alpha
>0,~p>2, 
\]
resting upon compatibility conditions between $\alpha $ and $p$.
Letting $2_{\alpha }:=2N/(N-\alpha)$ and
$2_{\alpha }^{*}:=2(2N-2+\alpha)/(2N-2-\alpha)$,
the problem is still open for $0<\alpha <2$ and $2_{\alpha
}<p\leq 2_{\alpha }^{*}$, for $2<\alpha <N$ and $2_{\alpha }^{*}\leq
p<2_{\alpha }$, and for $N\leq \alpha <2N-2$ and $p\geq 2_{\alpha }^{*}$.
Here we give a negative answer to the problem of the existence of radial
solutions in the first open case.
\end{abstract}


\maketitle

\section{Introduction}

In this paper we consider the following nonlinear problem: 
\begin{equation}
\left\{ 
\begin{array}{l}
\begin{array}{ll}
\smallskip -\triangle u+\displaystyle\dfrac{A}{\left| x\right| ^{\alpha }}u=u^{p-1}~ & 
\mathrm{in~}\mathbb{R}^{N}\setminus \left\{ 0\right\} ,~N\geq 3 \\ 
\medskip u>0 & \mathrm{in~}\mathbb{R}^{N}\setminus \left\{ 0\right\} 
\end{array}
\\ 
\begin{array}{l}
u\in H_{\alpha }^{1}\cap L^{p}(\mathbb{R}^{N})
\end{array}
\end{array}
\right.   \label{P}
\end{equation}
where $A,\alpha >0$, $p>2$ and $H_{\alpha }^{1}:=D^{1,2}(\mathbb{R}^{N})\cap
L^{2}(\mathbb{R}^{N},\left| x\right| ^{-\alpha }dx)$ is the natural energy
space related to the equation. We deal with problem (\ref{P})
in the classical sense, that is, speaking about \emph{solutions}
to (\ref{P}) we will always mean \emph{classical solutions%
} (cf. Remark \ref{RMK(classical)} below).

Problems like (\ref{P}) arise for instance in the search of solitary waves
for nonlinear Schr\"{o}dinger and Klein-Gordon equations with potential (see
e.g. \cite[Chapter 7]{THESIS}, \cite{BRnorm}, \cite{BRcharge}, the overviews in \cite{BBR 1}, \cite{BF sw
intro} and the monographs \cite{Strauss L.Notes}, \cite{YangY}) and (\ref{P}%
) itself is a radial model problem for the so-called \emph{zero mass} case
(see \cite{BPR}, \cite{BR TMA} and the references therein). In this respect, the
requirement $u\in H_{\alpha }^{1}\cap L^{p}(\mathbb{R}^{N})$ plays a
preeminent role, since it is necessary for the energy of the particle
represented by the solution to be finite.

Though it can be considered of quite recent investigation, problem (\ref{P})
has already some history and several existence and nonexistence results are
known, resting upon compatibility conditions between $\alpha $ and $p$ (see 
\cite{BGR} for a related cylindrical problem). At our knowledge, the first
results are due to Terracini \cite{Terracini}, who both proved that (\ref{P}%
) has no solution if 
\[
\left\{ 
\begin{array}{l}
\alpha =2 \\ 
p\neq 2^{*}
\end{array}
\right. \quad \mathrm{or}\quad \left\{ 
\begin{array}{l}
\alpha \neq 2 \\ 
p=2^{*}
\end{array}
\right. ,\qquad 2^{*}:=\frac{2N}{N-2},
\]
and explicitly found all the radial solutions of (\ref{P}) for $\left(
\alpha ,p\right) =\left( 2,2^{*}\right) $. As usual, $2^{*}$ denotes the
critical exponent for the Sobolev embedding in dimension $N\geq 3$. The
problem was subsequently addressed in \cite{Co-Cr-Par}, where it was proved
that (\ref{P}) has no solution if 
\[
\left\{ 
\begin{array}{l}
0<\alpha <2 \\ 
p>2^{*}
\end{array}
\right. \quad \mathrm{or}\quad \left\{ 
\begin{array}{l}
\alpha >2 \\ 
2<p<2^{*}
\end{array}
\right. .
\]
On the other hand, the authors obtained the existence of a radial solution
to (\ref{P}) provided that 
\[
\left\{ 
\begin{array}{l}
\smallskip 0<\alpha <2 \\ 
2^{*}+\frac{\alpha -2}{N-2}<p<2^{*}
\end{array}
\right. \quad \mathrm{or}\quad \left\{ 
\begin{array}{l}
\smallskip \alpha >2 \\ 
2^{*}<p<2^{*}+\frac{\alpha -2}{N-2}
\end{array}
\right. .
\]
The existence and nonexistence results of \cite{Co-Cr-Par} were then
extended in \cite{BRpow}, by showing that (\ref{P}) has no solution also if 
\[
\left\{ 
\begin{array}{l}
0<\alpha <2 \\ 
2<p\leq 2_{\alpha }
\end{array}
\right. \quad \mathrm{or}\quad \left\{ 
\begin{array}{l}
2<\alpha <N \\ 
p\geq 2_{\alpha }
\end{array}
\right. ,\qquad 2_{\alpha }:=\frac{2N}{N-\alpha },
\]
and obtaining a radial solution for every pair $\left( \alpha ,p\right) $
such that 
\[
\left\{ 
\begin{array}{l}
\smallskip 0<\alpha <2 \\ 
2^{*}+2\frac{\alpha -2}{N-2}<p<2^{*}
\end{array}
\right. \quad \mathrm{or}\quad \left\{ 
\begin{array}{l}
\smallskip \alpha >2 \\ 
2^{*}<p<2^{*}+2\frac{\alpha -2}{N-2}
\end{array}
\right. .
\]
A further extension of this existence condition were found in \cite
{Su-Wang-Will 2}, \cite{Su-Wang-Will p}, where the authors proved that (\ref
{P}) has a radial solution for all the pairs $\left( \alpha ,p\right) $
satisfying 
\[
\left\{ 
\begin{array}{l}
0<\alpha <2 \\ 
2_{\alpha }^{*}<p<2^{*}
\end{array}
\right. \quad \mathrm{or}\quad \left\{ 
\begin{array}{l}
2<\alpha <2N-2 \\ 
2^{*}<p<2_{\alpha }^{*}
\end{array}
\right. \quad \mathrm{or}\quad \left\{ 
\begin{array}{l}
\alpha \geq 2N-2 \\ 
p>2^{*}
\end{array}
\right. ,\quad 2_{\alpha }^{*}:=2\frac{2N-2+\alpha }{2N-2-\alpha }.
\]

All these known results are portrayed in the below picture of the $\alpha p$%
-plane, where the nonexistence regions are shaded in light gray (and include
both the lines $p=2^{*}$ and $p=2_{\alpha }$, except for the pair $\left(
\alpha ,p\right) =\left( 2,2^{*}\right) $), while dark gray means existence.
The problem is still open for the pairs $\left( \alpha ,p\right) $ in the
white regions of the picture, namely, for 
\[
\left\{ 
\begin{array}{l}
0<\alpha <2 \\ 
2_{\alpha }<p\leq 2_{\alpha }^{*}
\end{array}
\right. ,\quad \quad \left\{ 
\begin{array}{l}
2<\alpha <N \\ 
2_{\alpha }^{*}\leq p<2_{\alpha }
\end{array}
\right. \quad \mathrm{and}\quad \left\{ 
\begin{array}{l}
N\leq \alpha <2N-2 \\ 
p\geq 2_{\alpha }^{*}
\end{array}
\right. .
\]
In this paper we give a negative answer to the problem of radial solutions
to (\ref{P}) in the first of such cases.

\begin{figure}[ht]
   \centering
   \includegraphics{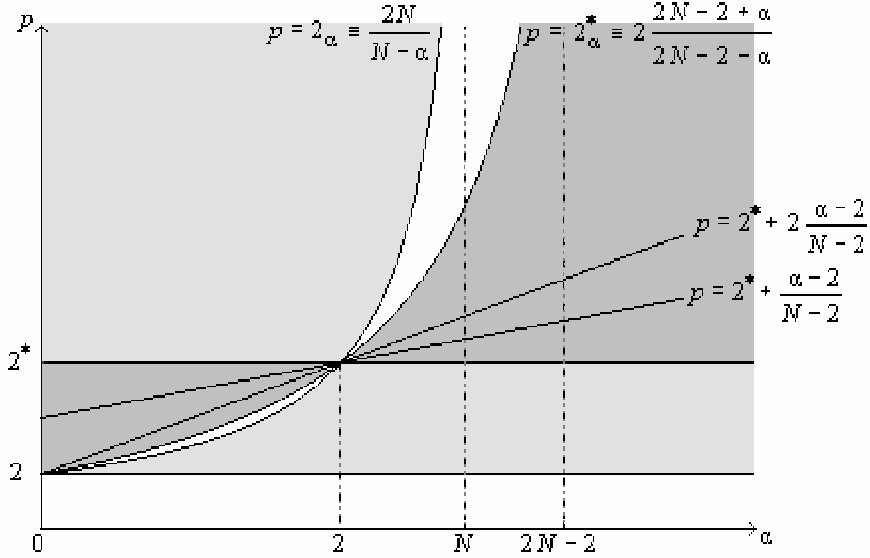}
   \caption{Existence (dark gray) and nonexistence (light gray) regions.}
\end{figure}

\begin{remark}
\label{RMK(classical)}In some of the above mentioned existence results, the
authors only concern themselves with nonnegative radial weak solutions to (%
\ref{P}) in the sense of $H_{\alpha }^{-1}$, the dual space of $H_{\alpha
}^{1}$ (in \cite{Su-Wang-Will 2}, \cite{Su-Wang-Will p}, only weak solutions
in the sense of the dual space of the radial subspace of $H_{\alpha }^{1}$
are considered, but the symmetric criticality type results of \cite{BGRrmks}
apply, yielding solutions in the sense of $H_{\alpha }^{-1}$). However,
Schauder regularity theory and the strong maximum principle (c.f. also
Remark \ref{COR(v>0)} below) assure that all such solutions are actually
positive classical solutions to (\ref{P}).
\end{remark}

Our nonexistence result is the following.

\begin{theorem}
\label{THM(nonex)}Let $0<\alpha <2$ and $2_{\alpha }<p\leq 2_{\alpha }^{*}$.
Then (\ref{P}) has no radial solution.
\end{theorem}

Observe that, although we are concerned with classical solutions, Theorem 
\ref{THM(nonex)} also prevents the existence of nonnegative radial weak
solutions in the $H_{\alpha }^{-1}$ sense, by the same reasons used in
Remark \ref{RMK(classical)}.

Theorem \ref{THM(nonex)} will be proved in Section \ref{SEC:nonex} and needs
a more refined argument than the one used in \cite{BRpow}, \cite{Co-Cr-Par}, 
\cite{Terracini}, where the nonexistence results were all obtained by
Pohozaev type identities. In fact, we will combine a Pucci-Serrin type
identity (see \cite{Pucci-Serrin}), which we deduce by an argument of \cite
{Co-Terr-Verz}, with a suitable asymptotic estimate (Lemma \ref{LEM(liminf)}%
), which derives from our next result.

\begin{theorem}
\label{THM(stime)}Let $0<\alpha <2$ and $2_{\alpha }<p<2^{*}$. Assume that $u
$ is a radial solution of 
\begin{equation}
\left\{ 
\begin{array}{l}
\begin{array}{ll}
\smallskip -\triangle u+\dfrac{A}{\left| x\right| ^{\alpha }}u=u^{p-1}~ & 
\mathrm{in~}\mathbb{R}^{N}\setminus \left\{ 0\right\} ,~N\geq 3 \\ 
\medskip u>0 & \mathrm{in~}\mathbb{R}^{N}\setminus \left\{ 0\right\} 
\end{array}
\\ 
\begin{array}{l}
u\in H_{\alpha }^{1}
\end{array}
\end{array}
\right. .  \label{P1}
\end{equation}
Then, as $x\rightarrow 0$, one has 
\begin{equation}
u\left( x\right) =\left\{ 
\begin{array}{ll}
\medskip O\left( 1\right)  & \mathrm{if~}p<2^{*}-1 \\ 
\medskip O\left( \ln \left| x\right| \right)  & \mathrm{if~}p=2^{*}-1 \\ 
O\left( \left| x\right| ^{-\frac{N-2}{2}(p-2^{*}+1)}\right) \quad  & \mathrm{%
if~}p>2^{*}-1
\end{array}
\right. .  \label{THM(stime):stime}
\end{equation}
\end{theorem}

Theorem \ref{THM(stime)} will be proved in Section \ref{SEC:stime} and,
besides yielding Theorem \ref{THM(nonex)}, it is interesting on its own,
since it also covers the existence case $2_{\alpha }^{*}<p<2^{*}$ (some
results on the asymptotic behaviour of solutions at infinity can be found in 
\cite{THESIS}, \cite{G-R}). Observe that all the cases of (\ref
{THM(stime):stime}) improve the estimate of a well known Radial Lemma for $%
D^{1,2}(\mathbb{R}^{N})$ (see \cite[Lemma A.III]{Beres-Lions}, where the
proof also works for $0<\left| x\right| <1$). Moreover, they are all
possible for $2_{\alpha }<p<2^{*}$ (even for $2_{\alpha }<p\leq 2_{\alpha
}^{*}$) if $N<6$, whereas only the third case occurs if $N\geq 7$.

Our proof of Theorem \ref{THM(stime)} will proceed as follows. First, we
will consider the ODE problem associated to the radial solutions of (\ref{P1}%
) and, after rescaling, we will recover its solutions as fixed points of a
suitable integral operator, which is expressed in terms of the modified
Bessel functions of the first and second kind (Lemma \ref{LEM(ptofisso)}).
Then we will show that such fixed points need to satisfy suitable estimates
(Theorem \ref{THM(stime_v)}), by exploiting a version of the already
mentioned Radial Lemma (Lemma \ref{LEM(pointwise)}), the monotonicity of the
integral operator and the well known behaviour of the Bessel functions at
the origin. Such estimates yield Theorem \ref{THM(stime)} by rescaling back.

Some useful properties of the modified Bessel functions are collected in the
Appendix. For a complete treatment, we refer the reader to \cite{Niki-Uva}, 
\cite{Temme} and the monumental monograph \cite{Watson}.
\medskip

\noindent \textbf{Notations. }We end this introductory section by
summarizing the notations of most frequent use throughout the
paper.

\noindent $\bullet $ We denote by $2^{*}:=2N/(N-2)$ the critical exponent
for the Sobolev embedding in dimension $N\geq 3$. Moreover we denote $%
2_{\alpha }:=2N/(N-\alpha )$ and $2_{\alpha }^{*}:=2(2N-2+\alpha
)/(2N-2-\alpha )$.

\noindent $\bullet $ We set $\mathbb{R}_{+}:=\left( 0,+\infty \right) .$

\noindent $\bullet $ If $\Omega \subseteq \mathbb{R}^{d}$, $d\geq 1$, is a
measurable set, $\rho :\Omega \rightarrow \mathbb{R}_{+}$ is a measurable
function and $1\leq q\leq \infty $, then $L^{q}(\Omega ,\rho \left( z\right)
dz)$ is the usual real Lebesgue space with respect to the measure $\rho
\left( z\right) dz$ ($dz$ stands for the Lebesgue measure on $\mathbb{R}^{d}$%
).

\noindent $\bullet $ $D^{1,2}(\mathbb{R}^{N})=\{u\in L^{2^{*}}(\mathbb{R}%
^{N}):\nabla u\in L^{2}(\mathbb{R}^{N})\}$ is the usual Sobolev space, which
identifies with the completion of $C_{\mathrm{c}}^{\infty }(\mathbb{R}^{N})$
with respect to the norm of the gradient.

\noindent $\bullet $ $I_{\nu }$ and $K_{\nu }$ are the modified Bessel
functions of order $\nu $, of the first and second kind
respectively.

\noindent $\bullet $ $o$ and $O$ are the usual Landau symbols. Moreover, by $%
f\left( t\right) \sim g\left( t\right) $ and $f\left( t\right) \asymp
g\left( t\right) $ as $t\rightarrow t_{0}$ we respectively mean $%
\displaystyle \lim_{t\rightarrow t_{0}}f\left( t\right) /g\left( t\right) =1$
and $\displaystyle \lim_{t\rightarrow t_{0}}f\left( t\right) /g\left(
t\right) =\ell \in \mathbb{R}\setminus \left\{ 0\right\} $.

\section{Asymptotic estimates for radial solutions at the origin\label%
{SEC:stime}}

In this section we assume $0<\alpha <2$ and $2_{\alpha }<p<2^{*}$. As one
can easily check, the problem of radial solutions $u\left( x\right) =\phi
\left( \left| x\right| \right) $ to (\ref{P1}) is equivalent to the
following ODE problem:

\begin{equation}
\left\{ 
\begin{array}{ll}
\smallskip -\phi ^{\prime \prime }-\dfrac{N-1}{r}\phi ^{\prime }+\dfrac{A}{%
r^{\alpha }}\phi =\phi ^{p-1} & \mathrm{in~}\mathbb{R}_{+}=\left( 0,+\infty
\right) \\ 
\medskip \phi >0 & \mathrm{in~}\mathbb{R}_{+} \\ 
r^{-\frac{\alpha }{2}}\phi ,\,\phi ^{\prime }\in L^{2}(\mathbb{R}%
_{+},r^{N-1}dr) & 
\end{array}
\right.  \label{Pb_phi}
\end{equation}
(cf. the proof of Lemma \ref{LEM(pointwise)} below). Making the change of
variable 
\begin{equation}
r=r\left( t\right) =\left( \frac{2-\alpha }{2\sqrt{A}}t\right) ^{\frac{2}{%
2-\alpha }}  \label{cambio1}
\end{equation}
and defining 
\begin{equation}
v\left( t\right) =\phi \left( r\left( t\right) \right) \quad \mathrm{for~all~%
}t>0,  \label{cambio2}
\end{equation}
one has 
\begin{equation}
t=t\left( r\right) =\frac{2\sqrt{A}}{2-\alpha }r^{\frac{2-\alpha }{2}%
},\qquad \phi \left( r\right) =v\left( t\left( r\right) \right) =v\left( 
\frac{2\sqrt{A}}{2-\alpha }r^{\frac{2-\alpha }{2}}\right) ,  \label{cambio3}
\end{equation}
so that 
\[
\phi ^{\prime }\left( r\right) =v^{\prime }\left( t\right) \frac{dt}{dr}=%
\sqrt{A}v^{\prime }\left( t\right) r^{-\frac{\alpha }{2}} 
\]
and 
\begin{eqnarray*}
\phi ^{\prime \prime }\left( r\right) &=&\sqrt{A}\left( v^{\prime \prime
}\left( t\right) \frac{dt}{dr}r^{-\frac{\alpha }{2}}-\frac{\alpha }{2}%
v^{\prime }\left( t\right) r^{-\frac{\alpha +2}{2}}\right) =\sqrt{A}\left(
v^{\prime \prime }\left( t\right) \sqrt{A}r^{-\alpha }-\frac{\alpha }{2}%
v^{\prime }\left( t\right) r^{-\frac{\alpha +2}{2}}\right) \\
&=&Av^{\prime \prime }\left( t\right) r^{-\alpha }-\sqrt{A}\frac{\alpha }{2}%
v^{\prime }\left( t\right) r^{-\frac{\alpha +2}{2}}.
\end{eqnarray*}
Plugging into the equation of (\ref{Pb_phi}) we get 
\[
-Av^{\prime \prime }\left( t\right) r^{-\alpha }+\left( \frac{\alpha }{2}%
-N+1\right) \sqrt{A}v^{\prime }\left( t\right) r^{-\frac{\alpha +2}{2}}+%
\dfrac{A}{r^{\alpha }}v\left( t\right) =v\left( t\right) ^{p-1} 
\]
and multiplying both sides by $r^{\alpha }/A$ we obtain 
\[
-v^{\prime \prime }\left( t\right) +\frac{\alpha -2N+2}{2\sqrt{A}}v^{\prime
}\left( t\right) r^{\frac{\alpha -2}{2}}+v\left( t\right) =\frac{r^{\alpha }%
}{A}v\left( t\right) ^{p-1}. 
\]
Since $r^{\alpha }=\left( \frac{2-\alpha }{2\sqrt{A}}\right) ^{\frac{2\alpha 
}{2-\alpha }}t^{\frac{2\alpha }{2-\alpha }}$ and $r^{\frac{\alpha -2}{2}}=%
\frac{2\sqrt{A}}{2-\alpha }\frac{1}{t}$, the equation of (\ref{Pb_phi})
turns thus out to be equivalent to 
\[
-v^{\prime \prime }-\frac{2N-2-\alpha }{2-\alpha }\frac{1}{t}v^{\prime
}+v=\left( \frac{2-\alpha }{2A^{1/\alpha }}\right) ^{\frac{2\alpha }{%
2-\alpha }}t^{\frac{2\alpha }{2-\alpha }}v^{p-1}. 
\]
Observing that 
\[
r^{N-1-\alpha }dr=\left( \frac{2-\alpha }{2\sqrt{A}}t\right) ^{\frac{2\left(
N-1-\alpha \right) }{2-\alpha }}\frac{1}{\sqrt{A}}\left( \frac{2-\alpha }{2%
\sqrt{A}}t\right) ^{\frac{\alpha }{2-\alpha }}dt=\left( \mathrm{const.}%
\right) t^{\frac{2N-2-\alpha }{2-\alpha }}dt 
\]
and $\phi ^{\prime }\left( r\right) ^{2}=Av^{\prime }\left( t\right)
^{2}r^{-\alpha }$, one has 
\[
\int_{0}^{+\infty }\phi \left( r\right) ^{2}r^{N-1-\alpha }dr=\left( \mathrm{%
const.}\right) \int_{0}^{+\infty }v\left( t\right) ^{2}t^{\frac{2N-2-\alpha 
}{2-\alpha }}dt 
\]
and 
\begin{equation}
\int_{0}^{+\infty }\phi ^{\prime }\left( r\right)
^{2}r^{N-1}dr=A\int_{0}^{+\infty }v^{\prime }\left( t\left( r\right) \right)
^{2}r^{N-1-\alpha }dr=\left( \mathrm{const.}\right) \int_{0}^{+\infty
}v^{\prime }\left( t\right) ^{2}t^{\frac{2N-2-\alpha }{2-\alpha }}dt.
\label{int_phi'}
\end{equation}
As a conclusion, setting 
\[
\nu :=\frac{N-2}{2-\alpha },\qquad B:=\left( \frac{2-\alpha }{2A^{1/\alpha }}%
\right) ^{\frac{2\alpha }{2-\alpha }}, 
\]
problem (\ref{Pb_phi}) is equivalent to 
\begin{equation}
\left\{ 
\begin{array}{ll}
\smallskip -v^{\prime \prime }-\dfrac{2\nu +1}{t}v^{\prime }+v=Bt^{\frac{%
2\alpha }{2-\alpha }}v^{p-1} & \mathrm{in~}\mathbb{R}_{+} \\ 
\smallskip v>0 & \mathrm{in~}\mathbb{R}_{+} \\ 
v\in H & 
\end{array}
\right.  \label{Pb_v}
\end{equation}
where 
\[
H:=H^{1}(\mathbb{R}_{+},t^{\frac{2N-2-\alpha }{2-\alpha }}dt):=\left\{ v\in
L^{2}(\mathbb{R}_{+},t^{\frac{2N-2-\alpha }{2-\alpha }}dt):v^{\prime }\in
L^{2}(\mathbb{R}_{+},t^{\frac{2N-2-\alpha }{2-\alpha }}dt)\right\} . 
\]
Note that $\nu ,B>0$.

The next lemma is a version of a well known Radial Lemma \cite{Beres-Lions}
and states some properties of the functions in $H$.

\begin{lemma}
\label{LEM(pointwise)}Every $v\in H$ is continuous on $\mathbb{R}_{+}$ (up
to the choice of a representative) and satisfies 
\begin{equation}
\left| v\left( t\right) \right| \leq C_{N,A,\alpha }\left\| v^{\prime
}\right\| _{2,\alpha }\frac{1}{t^{\nu }}\quad \mathrm{for~all~}t>0,
\label{LEM(pointwise): stima}
\end{equation}
where $\left\| v^{\prime }\right\| _{2,\alpha }$ is the norm of $v^{\prime }$
in $L^{2}(\mathbb{R}_{+},t^{\frac{2N-2-\alpha }{2-\alpha }}dt)$ and the
constant $C_{N,A,\alpha }$ only depends on $N$, $A$ and $\alpha $.
\end{lemma}

\proof%
Let $v\in H$ and let $\phi $ be defined by (\ref{cambio1})-(\ref{cambio3}).
Then $\phi \left( \left| x\right| \right) $ belongs to $D^{1,2}(\mathbb{R}%
^{N})=\{u\in L^{2^{*}}(\mathbb{R}^{N}):\nabla u\in L^{2}(\mathbb{R}^{N})\}$.
Indeed $r^{-\frac{\alpha }{2}}\phi \in L^{2}(\mathbb{R}_{+},r^{N-1}dr)$
implies that $\phi \left( \left| x\right| \right) \in L^{2}(\mathbb{R}%
^{N},\left| x\right| ^{-\alpha }dx)\subset L_{\mathrm{loc}}^{1}(\mathbb{R}%
^{N})$ and $\phi \in L^{2}(\left( 1,+\infty \right) )$, while $\phi ^{\prime
}\in L^{2}(\mathbb{R}_{+},r^{N-1}dr)$ implies $\phi ^{\prime }\in
L^{2}(\left( 1,+\infty \right) )$, as well as that the gradient of $\phi
\left( \left| x\right| \right) $ is in $L^{2}(\mathbb{R}^{N})$; hence $\phi
\in H^{1}(\left( 1,+\infty \right) )$ and thus $\lim_{\left| x\right|
\rightarrow \infty }\phi \left( \left| x\right| \right) =0$, which yields $%
u\in L^{2^{*}}(\mathbb{R}^{N})$ by Sobolev inequality (see the version given
in \cite[Theorem 8.3]{Lieb-Loss}).

So, by \cite[Lemma A.III]{Beres-Lions} (where the proof actually works for
every $x\neq 0$), $\phi $ is continuous on $\mathbb{R}_{+}$ (up to the
choice of a representative) and satisfies 
\[
\left| \phi \left( r\right) \right| \leq C_{N}\left( \int_{0}^{+\infty }\phi
^{\prime }\left( r\right) ^{2}r^{N-1}dr\right) ^{1/2}\frac{1}{r^{\frac{N-2}{2%
}}}\quad \mathrm{for~all~}r>0, 
\]
where the constant $C_{N}$ only depends on $N$. This gives (\ref
{LEM(pointwise): stima}) by (\ref{int_phi'}) and (\ref{cambio1}).%
\endproof%
\bigskip

We now consider the linear equation associated to the equation of (\ref{Pb_v}
), whose general solution can be expressed in terms of the modified Bessel
functions of the first and second kind (see Appendix).

\begin{lemma}
\label{LEM(IntGen)}For any $g\in C\left( \mathbb{R}_{+}\right) $, the
general solution of the equation 
\begin{equation}
-v^{\prime \prime }-\dfrac{2\nu +1}{t}v^{\prime }+v=g\left( t\right) \quad 
\mathrm{in~}\mathbb{R}_{+}  \label{eq_lin}
\end{equation}
is 
\begin{eqnarray}
&&v\left( t;c_{1},c_{2}\right)=   \label{IntGen} \\
&=&\frac{1}{t^{\nu }}\left\{ \left( c_{1}-\int_{1}^{t}s^{1+\nu }K_{\nu
}\left( s\right) g\left( s\right) ds\right) I_{\nu }\left( t\right) +\left(
c_{2}+\int_{1}^{t}s^{1+\nu }I_{\nu }\left( s\right) g\left( s\right)
ds\right) K_{\nu }\left( t\right) \right\} ,  \nonumber
\end{eqnarray}
where $c_{1},c_{2}\in \mathbb{R}$ are arbitrary constants and $I_{\nu }$ and 
$K_{\nu }$ are the modified Bessel functions of order $\nu $, of the first
and second kind respectively.
\end{lemma}

\proof%
Taking into account that $I_{\nu }$ and $K_{\nu }$ are linearly independent
solutions of the modified Bessel equation 
\[
-v^{\prime \prime }-\dfrac{1}{t}v^{\prime }+\left( 1+\dfrac{\nu ^{2}}{t^{2}}%
\right) v=0\quad \mathrm{in~}\mathbb{R}_{+}, 
\]
one easily checks that $t^{-\nu }I_{\nu }$ and $t^{-\nu }K_{\nu }$ are
linearly independent solutions of the homogeneous equation associated to (%
\ref{eq_lin}). On the other hand, the function 
\[
\widetilde{v}\left( t\right) =\frac{1}{t^{\nu }}\left( K_{\nu }\left(
t\right) \int_{1}^{t}s^{1+\nu }I_{\nu }\left( s\right) g\left( s\right)
ds-I_{\nu }\left( t\right) \int_{1}^{t}s^{1+\nu }K_{\nu }\left( s\right)
g\left( s\right) ds\right) 
\]
is a particular solution of equation (\ref{eq_lin}), since one has 
\begin{eqnarray*}
&&-\widetilde{v}^{\prime \prime }\left( t\right) -\dfrac{2\nu +1}{t}%
\widetilde{v}^{\prime }\left( t\right) +\widetilde{v}\left( t\right)  \\
&=&-\frac{1}{t^{\nu }}\left( K_{\nu }^{\prime \prime }\left( t\right) +\frac{%
1}{t}K_{\nu }^{\prime }\left( t\right) -\left( 1+\frac{\nu ^{2}}{t^{2}}%
\right) K_{\nu }\left( t\right) \right) \int_{1}^{t}s^{1+\nu }I_{\nu }\left(
s\right) g\left( s\right) ds+ \\
&&+\frac{1}{t^{\nu }}\left( I_{\nu }^{\prime \prime }\left( t\right) +\frac{1%
}{t}I_{\nu }^{\prime }\left( t\right) -\left( 1+\frac{\nu ^{2}}{t^{2}}%
\right) I_{\nu }\left( t\right) \right) \int_{1}^{t}s^{1+\nu }K_{\nu }\left(
s\right) g\left( s\right) ds+ \\
&&+t\left( I_{\nu }\left( t\right) K_{\nu +1}\left( t\right) +K_{\nu }\left(
t\right) I_{\nu +1}\left( t\right) \right) g\left( t\right) 
\end{eqnarray*}
and the following identity holds: $I_{\nu }\left( t\right) K_{\nu +1}\left(
t\right) +K_{\nu }\left( t\right) I_{\nu +1}\left( t\right) =\frac{1}{t}$
for all $t>0$.%
\endproof%
\bigskip

In the following, for the sake of brevity, we will denote 
\[
H_{+}:=\left\{ v\in H:v>0\right\} 
\]
and 
\[
I\left( t\right) :=t^{\frac{N+\alpha }{2-\alpha }}I_{\nu }\left( t\right)
\quad \mathrm{and}\quad K\left( t\right) :=t^{\frac{N+\alpha }{2-\alpha }%
}K_{\nu }\left( t\right) \quad \mathrm{for~every~}t>0. 
\]
Furthermore, we will make an extensive use of the following estimates (see
the Appendix for more accurate asymptotic equivalences):
\medskip 

\noindent $\bullet$ as $t\rightarrow 0^{+}$ one has
\begin{equation}
\frac{I_{\nu +1}\left( t\right) }{t}\asymp I_{\nu }\left( t\right) \asymp t^{\nu },\ \ 
tK_{\nu +1}\left( t\right) \asymp K_{\nu }\left( t\right) \asymp t^{-\nu },\ \ 
I\left( t\right) \asymp t^{\frac{N+\alpha }{2-\alpha }+\nu},\ \ 
K\left( t\right) \asymp t^{\frac{N+\alpha }{2-\alpha }-\nu };
\label{IKzero}
\end{equation}

\noindent $\bullet$ as $t\rightarrow +\infty $ one has
\begin{equation}
I_{\nu }\left( t\right) \asymp \frac{e^{t}}{\sqrt{t}},\ \ 
K_{\nu }\left(t\right) \asymp \frac{e^{-t}}{\sqrt{t}},\ \ 
I\left( t\right) \asymp t^{\frac{N+\alpha }{2-\alpha }-\frac{1}{2}}e^{t},\ \ 
K\left( t\right) \asymp t^{\frac{N+\alpha }{2-\alpha }-\frac{1}{2}}e^{-t}.  
\label{IKinfty}
\end{equation}

\noindent Note that $(N+\alpha )/(2-\alpha )=\nu +1+2\alpha /(2-\alpha )$.

\begin{lemma}
\label{LEM(ptofisso)}Let $v\in H_{+}$. Then $v$ is a solution to problem (%
\ref{Pb_v}) if and only if 
\begin{equation}
v\left( t\right) =\frac{B}{t^{\nu }}\left\{ I_{\nu }\left( t\right)
\int_{t}^{+\infty }K\left( s\right) v\left( s\right) ^{p-1}ds+K_{\nu }\left(
t\right) \int_{0}^{t}I\left( s\right) v\left( s\right) ^{p-1}ds\right\}
\quad \mathrm{for~all~}t>0.  \label{v=Tv}
\end{equation}
\end{lemma}

\begin{remark}
\label{RMK(T definito)}The integrals $\int_{t}^{+\infty }K\left( s\right)
v\left( s\right) ^{p-1}ds$ and $\int_{0}^{t}I\left( s\right) v\left(
s\right) ^{p-1}ds$ are finite for every $v\in H_{+}$ and $t>0$, since:

\begin{itemize}
\item  $K\left( s\right) \asymp s^{\frac{N+\alpha }{2-\alpha }-\frac{1}{2}%
}e^{-s}$ and $v\left( s\right) =O\left( s^{-\nu }\right) $ as $s\rightarrow
+\infty $ (see (\ref{IKinfty}) and Lemma \ref{LEM(pointwise)});

\item  from (\ref{IKzero}) and Lemma \ref{LEM(pointwise)}, it follows that 
\[
I\left( s\right) v\left( s\right) ^{p-1}\asymp s^{\frac{N+\alpha }{2-\alpha }%
+\nu }v\left( s\right) ^{p-1}=s^{\frac{N+\alpha }{2-\alpha }+\nu }O\left(
s^{-\nu (p-1)}\right) =O\left( s^{\frac{N+\alpha }{2-\alpha }+2\nu -\nu
p}\right) 
\]
as $s\rightarrow 0^+$, where 
\[
\frac{N+\alpha }{2-\alpha }+2\nu -\nu p+1=\nu \left( 2^{*}+1-p\right) >0.
\]
\end{itemize}
\end{remark}

\proof%
Clearly, $v$ solves (\ref{Pb_v}) if (\ref{v=Tv}) holds, since for all $t>0$
one has 
\begin{eqnarray*}
v\left( t\right) &=&\frac{I_{\nu }\left( t\right) }{t^{\nu }}\left(
B\int_{1}^{+\infty }K\left( s\right) v\left( s\right)
^{p-1}ds-B\int_{1}^{t}K\left( s\right) v\left( s\right) ^{p-1}ds\right) + \\
&&+\frac{K_{\nu }\left( t\right) }{t^{\nu }}\left( B\int_{0}^{1}I\left(
s\right) v\left( s\right) ^{p-1}ds+B\int_{1}^{t}I\left( s\right) v\left(
s\right) ^{p-1}ds\right)
\end{eqnarray*}
and thus $v$ is of the form (\ref{IntGen}) with $g\left( t\right)
=Bt^{2\alpha /(2-\alpha )}v\left( t\right) ^{p-1}$ continuous on $\mathbb{R}%
_{+}$. In order to prove the ``only if'' part of the lemma, assume that $v$
is a solution of problem (\ref{Pb_v}). Then, using Lemma \ref{LEM(IntGen)}
with $g\left( t\right) =Bt^{2\alpha /(2-\alpha )}v\left( t\right) ^{p-1}$,
there exist two unique constants $c_{1}=c_{1}\left( v\right)
,c_{2}=c_{2}\left( v\right) \in \mathbb{R}$ such that 
\begin{equation}
v\left( t\right) =\frac{1}{t^{\nu }}\left\{ \left(
c_{1}-B\int_{1}^{t}K\left( s\right) v^{p-1}\left( s\right) ds\right) I_{\nu
}\left( t\right) +\left( c_{2}+B\int_{1}^{t}I\left( s\right) v^{p-1}\left(
s\right) ds\right) K_{\nu }\left( t\right) \right\}
\label{v(t)=}
\end{equation}
for all $t>0$. Set 
\begin{eqnarray*}
\Phi _{1} &=&\Phi _{1}\left( t\right) :=\left( c_{1}-B\int_{1}^{t}K\left(
s\right) v^{p-1}\left( s\right) ds\right) \frac{I_{\nu }\left( t\right) }{%
t^{\nu }}, \\
\Phi _{2} &=&\Phi _{2}\left( t\right) :=\left( c_{2}+B\int_{1}^{t}I\left(
s\right) v^{p-1}\left( s\right) ds\right) \frac{K_{\nu }\left( t\right) }{%
t^{\nu }}
\end{eqnarray*}
in such a way that $v=\Phi _{1}+\Phi _{2}$, and assume by contradiction that 
\[
c_{1}\neq B_{1}:=B\int_{1}^{+\infty }K\left( s\right) v\left( s\right)
^{p-1}ds, 
\]
where $B_{1}<+\infty $ by Remark \ref{RMK(T definito)}. This implies 
\[
c_{1}-B\int_{1}^{t}K\left( s\right) v\left( s\right) ^{p-1}ds\asymp 1\quad 
\textrm{as }t\rightarrow +\infty 
\]
and hence, as $t\rightarrow +\infty $, one gets 
\begin{eqnarray*}
\Phi _{1}\Phi _{2} &=&t^{-2\nu }I_{\nu }\left( t\right) K_{\nu }\left(
t\right) \left( c_{1}-B\int_{1}^{t}K\left( s\right) v\left( s\right)
^{p-1}ds\right) \left( c_{2}+B\int_{1}^{t}I\left( s\right) v\left( s\right)
^{p-1}ds\right) \\
&\asymp &\left( c_{2}+B\int_{1}^{t}I\left( s\right) v\left( s\right)
^{p-1}ds\right) t^{-2\nu -1}
\end{eqnarray*}
and 
\[
\Phi _{1}^{2}=t^{-2\nu }I_{\nu }\left( t\right) ^{2}\left(
c_{1}-B\int_{1}^{t}K\left( s\right) v\left( s\right) ^{p-1}ds\right)
^{2}\asymp t^{-2\nu -1}e^{2t}. 
\]
Now we distinguish two cases, according to the value of the limit 
\[
\lim_{t\rightarrow +\infty }\left( c_{2}+B\int_{1}^{t}I\left(
s\right) v\left( s\right) ^{p-1}ds\right), 
\]
which exists since $I\left(
s\right) v\left( s\right) ^{p-1}>0$. If the limit is finite, we readily get 
\[
\Phi _{1}\Phi _{2}\asymp t^{-2\nu -1}\asymp e^{-2t}\Phi _{1}^{2}=o\left(
\Phi _{1}^{2}\right) \quad \mathrm{as~}t\rightarrow +\infty . 
\]
If the limit is infinite, then, by De L'H\^{o}pital's rule, we obtain 
\begin{eqnarray*}
\lim_{t\rightarrow +\infty }\frac{\Phi _{1}\Phi _{2}}{\Phi _{1}^{2}}
&=&\left( \mathrm{const.}\right) \lim_{t\rightarrow +\infty }\frac{%
c_{2}+B\int_{1}^{t}I\left( s\right) v\left( s\right) ^{p-1}ds}{e^{2t}}%
\stackrel{H}{=}\left( \mathrm{const.}\right) \lim_{t\rightarrow +\infty }%
\frac{I\left( t\right) v\left( t\right) ^{p-1}}{e^{2t}} \\
&=&\left( \mathrm{const.}\right) \lim_{t\rightarrow +\infty }\frac{t^{\frac{%
N+\alpha }{2-\alpha }-\frac{1}{2}}e^{t}O\left( t^{-\nu (p-1)}\right) }{e^{2t}%
}=0.
\end{eqnarray*}
So, in any case, we have $\Phi _{1}\Phi _{2}=o\left( \Phi _{1}^{2}\right) $
as $t\rightarrow +\infty $ and hence 
\[
v^{2}=\Phi _{1}^{2}+2\Phi _{1}\Phi _{2}+\Phi _{2}^{2}\geq \Phi
_{1}^{2}+2\Phi _{1}\Phi _{2}\sim \Phi _{1}^{2}\asymp t^{-2\nu -1}e^{2t}\quad 
\mathrm{as~}t\rightarrow +\infty . 
\]
This implies $v\notin L^{2}(\mathbb{R}_{+},t^{\frac{2N-2-\alpha }{2-\alpha }%
}dt)$, which is false by hypothesis, and thus it must be $c_{1}=B_{1}$.
Substituting into (\ref{v(t)=}), we obtain 
\begin{equation}
v\left( t\right) =\frac{1}{t^{\nu }}\left\{ BI_{\nu }\left( t\right)
\int_{t}^{+\infty }K\left( s\right) v\left( s\right) ^{p-1}ds+\left(
c_{2}+B\int_{1}^{t}I\left( s\right) v\left( s\right) ^{p-1}ds\right) K_{\nu
}\left( t\right) \right\}
\label{v(t)= bis}
\end{equation}
for all $t>0$. We now prove that 
\[
c_{2}=B_{2}:=B\int_{0}^{1}I\left( s\right) v\left( s\right) ^{p-1}ds, 
\]
where $B_{2}<+\infty $ by Remark \ref{RMK(T definito)}. Taking the
derivative of (\ref{v(t)= bis}) and using the identities 
\[
K_{\nu }^{\prime }\left( t\right) -\frac{\nu }{t}K_{\nu }\left( t\right)
=-K_{\nu +1}\left( t\right) ,\qquad I_{\nu }^{\prime }\left( t\right) -\frac{%
\nu }{t}I_{\nu }\left( t\right) =I_{\nu +1}\left( t\right) 
\]
and $I\left( t\right) K_{\nu }\left( t\right) -K\left( t\right) I_{\nu
}\left( t\right) =t^{\frac{N+\alpha }{2-\alpha }}I_{\nu }\left( t\right)
K_{\nu }\left( t\right) -t^{\frac{N+\alpha }{2-\alpha }}K_{\nu }\left(
t\right) I_{\nu }\left( t\right) =0$ on $\mathbb{R}_{+}$, we get 
\begin{eqnarray*}
&&v^{\prime }\left( t\right)= \\
&=&-\frac{\nu }{t^{\nu +1}}\left\{ BI_{\nu
}\left( t\right) \int_{t}^{+\infty }K\left( s\right) v\left( s\right)
^{p-1}ds+\left( c_{2}+B\int_{1}^{t}I\left( s\right) v\left( s\right)
^{p-1}ds\right) K_{\nu }\left( t\right) \right\} + \\
&&+\frac{1}{t^{\nu }}\left\{ BI_{\nu }^{\prime }\left( t\right)
\int_{t}^{+\infty }K\left( s\right) v\left( s\right) ^{p-1}ds-BI_{\nu
}\left( t\right) K\left( t\right) v\left( t\right) ^{p-1}\right\} + \\
&&+\frac{1}{t^{\nu }}\left\{ \left( c_{2}+B\int_{1}^{t}I\left( s\right)
v\left( s\right) ^{p-1}ds\right) K_{\nu }^{\prime }\left( t\right) +BK_{\nu
}\left( t\right) I\left( t\right) v\left( t\right) ^{p-1}\right\} \\
&=&\frac{B}{t^{\nu }}v\left( t\right) ^{p-1}\left( K_{\nu }\left( t\right)
I\left( t\right) -I_{\nu }\left( t\right) K\left( t\right) \right) +\\
&&+\frac{B}{t^{\nu }}\left( I_{\nu }^{\prime }\left( t\right) -\frac{\nu }{t}I_{\nu
}\left( t\right) \right) \int_{t}^{+\infty }K\left( s\right) v\left(
s\right) ^{p-1}ds+ \\
&&+\frac{1}{t^{\nu }}\left( K_{\nu }^{\prime }\left( t\right) -\frac{\nu }{t}%
K_{\nu }\left( t\right) \right) \left( c_{2}+B\int_{1}^{t}I\left( s\right)
v\left( s\right) ^{p-1}ds\right) \\
&=&\frac{B}{t^{\nu }}I_{\nu +1}\left( t\right) \int_{t}^{+\infty }K\left(
s\right) v\left( s\right) ^{p-1}ds-\frac{1}{t^{\nu }}K_{\nu +1}\left(
t\right) \left( c_{2}+B\int_{1}^{t}I\left( s\right) v\left( s\right)
^{p-1}ds\right) .
\end{eqnarray*}
Setting 
\begin{eqnarray}
\Psi _{1} &=&\Psi _{1}\left( t\right) :=\frac{I_{\nu +1}\left( t\right) }{%
t^{\nu }}\int_{t}^{+\infty }K\left( s\right) v\left( s\right) ^{p-1}ds,
\label{psi1:=} \\
\Psi _{2} &=&\Psi _{2}\left( t\right) :=\frac{K_{\nu +1}\left( t\right) }{%
t^{\nu }}\left( c_{2}+B\int_{1}^{t}I\left( s\right) v\left( s\right)
^{p-1}ds\right) ,  \nonumber
\end{eqnarray}
in such a way that $v^{\prime }=\Psi _{1}+\Psi _{2}$, we show that 
\begin{equation}
\Psi _{1}\in L^{2}\left( \left( 0,1\right) ,t^{\frac{2N-2-\alpha }{2-\alpha }%
}dt\right) .  \label{psi1_inL2}
\end{equation}
If $\int_{0}^{+\infty }K\left( s\right) v\left( s\right) ^{p-1}ds<+\infty $,
then one has 
\begin{equation}
\Psi _{1}\asymp \frac{I_{\nu +1}\left( t\right) }{t^{\nu }}\asymp t\quad 
\mathrm{as~}t\rightarrow 0^{+}  \label{psi1_equigrande}
\end{equation}
and hence 
\[
t^{\frac{2N-2-\alpha }{2-\alpha }}\Psi _{1}^{2}\asymp t^{\frac{2N-2-\alpha }{%
2-\alpha }+2}=t^{\frac{2N+2-3\alpha }{2-\alpha }}\quad \mathrm{as~}%
t\rightarrow 0^{+} 
\]
with 
\[
\frac{2N+2-3\alpha }{2-\alpha }+1=2\frac{N+2-2\alpha }{2-\alpha }>2\frac{N-2%
}{2-\alpha }>0, 
\]
which implies (\ref{psi1_inL2}). Otherwise, if $\int_{0}^{+\infty }K\left(
s\right) v\left( s\right) ^{p-1}ds=+\infty $, we observe that 
\[
\frac{\nu }{2}\left( \frac{2\alpha }{N-2}-p\right) <\frac{\nu }{2}\left( 
\frac{2\alpha }{N-2}-\frac{2N}{N-\alpha }\right) <0 
\]
and apply De L'H\^{o}pital's rule: we obtain 
\begin{eqnarray*}
\lim_{t\rightarrow 0^{+}}\frac{\int_{t}^{+\infty }K\left( s\right) v\left(
s\right) ^{p-1}ds}{t^{\frac{\nu }{2}\left( \frac{2\alpha }{N-2}-p\right) }}%
&\stackrel{H}{=}&\left( \mathrm{const.}\right) \lim_{t\rightarrow 0^{+}}\frac{%
K\left( t\right) v\left( t\right) ^{p-1}}{t^{\frac{\nu }{2}\left( \frac{%
2\alpha }{N-2}-p\right) -1}} \\
&=&\left( \mathrm{const.}\right)
\lim_{t\rightarrow 0^{+}}\frac{t^{\frac{N+\alpha }{2-\alpha }-\nu }O\left(
t^{-\nu (p-1)}\right) }{t^{\frac{\nu }{2}\left( \frac{2\alpha }{N-2}%
-p\right) -1}} \\
&=&\left( \mathrm{const.}\right) \lim_{t\rightarrow 0^{+}}O\left( t^{\frac{%
N+2-\alpha }{2-\alpha }-\frac{\nu }{2}p}\right) =0,
\end{eqnarray*}
since 
\[
\frac{N+2-\alpha }{2-\alpha }-\frac{\nu }{2}p=\frac{\nu }{2}\left( 2^{*}-p+2%
\frac{2-\alpha }{N-2}\right) >0. 
\]
So, recalling (\ref{psi1:=}) and (\ref{psi1_equigrande}), one has 
\[
\Psi _{1}=\frac{I_{\nu +1}\left( t\right) }{t^{\nu }}o\left( t^{\frac{\nu }{2%
}\left( \frac{2\alpha }{N-2}-p\right) }\right) =o\left( t^{\frac{\nu }{2}%
\left( \frac{2\alpha }{N-2}-p\right) +1}\right) \quad \mathrm{as~}%
t\rightarrow 0^{+} 
\]
and hence 
\[
t^{\frac{2N-2-\alpha }{2-\alpha }}\Psi _{1}^{2}=o\left( t^{\frac{2N-2-\alpha 
}{2-\alpha }+\nu \left( \frac{2\alpha }{N-2}-p\right) +2}\right) \quad 
\mathrm{as~}t\rightarrow 0^{+} 
\]
with 
\[
\frac{2N-2-\alpha }{2-\alpha }+\nu \left( \frac{2\alpha }{N-2}-p\right)
+3=\nu \left( 2^{*}-p+2\frac{2-\alpha }{N-2}\right) >0, 
\]
which gives (\ref{psi1_inL2}) again. Therefore $v^{\prime }\in L^{2}(\left(
0,1\right) ,t^{\frac{2N-2-\alpha }{2-\alpha }}dt)$ implies 
\[
\Psi _{2}\in
L^{2}\left(\left( 0,1\right) ,t^{\frac{2N-2-\alpha }{2-\alpha }}dt\right).
\]
But this is impossible if $c_{2}\neq B_{2}$, since $c_{2}\neq B_{2}$ implies 
\[
\Psi _{2}\sim \left( c_{2}-B_{2}\right) \frac{K_{\nu +1}\left( t\right) }{%
t^{\nu }}\asymp \frac{1}{t^{2\nu +1}}\quad \mathrm{as~}t\rightarrow 0^{+}, 
\]
whence 
\[
t^{\frac{2N-2-\alpha }{2-\alpha }}\Psi _{2}^{2}\asymp t^{\frac{2N-2-\alpha }{%
2-\alpha }-4\nu -2}\quad \mathrm{as~}t\rightarrow 0^{+} 
\]
with 
\[
\frac{2N-2-\alpha }{2-\alpha }-4\nu -1=-2\nu <0. 
\]
So it must be $c_{2}=B_{2}$ and (\ref{v=Tv}) then follows from (\ref{v(t)=
bis}).%
\endproof%

\begin{remark}
\label{COR(v>0)}Checking the proof of Lemma \ref{LEM(ptofisso)}, one readily
sees that (\ref{v=Tv}) also holds for every nonnegative $v\in H$ satisfying
equation (\ref{Pb_v}). This directly yields, without the use of the maximum
principle, that every nontrivial nonnegative solution $v\in H$ of equation (%
\ref{Pb_v}) is strictly positive on $\mathbb{R}_{+}$. Indeed, since $I_{\nu
}\left( t\right) ,K_{\nu }\left( t\right) ,I\left( t\right) ,K\left(
t\right) >0$ for all $t>0$, if there exists $t_{0}>0$ such that $v\left(
t_{0}\right) =0$ then (\ref{v=Tv}) implies 
\[
\int_{t_{0}}^{+\infty }K\left( s\right) v\left( s\right)
^{p-1}ds=\int_{0}^{t_{0}}I\left( s\right) v\left( s\right) ^{p-1}ds=0,
\]
which means $v=0$ on $\mathbb{R}_{+}$.
\end{remark}

\begin{theorem}
\label{THM(stime_v)}Assume that $v$ is a solution of problem (\ref{Pb_v}).
Then, as $t\rightarrow 0^{+}$, one has 
\begin{equation}
v\left( t\right) =\left\{ 
\begin{array}{ll}
\medskip O\left( 1\right)  & \mathrm{if~}p<2^{*}-1 \\ 
\medskip O\left( \ln t\right)  & \mathrm{if~}p=2^{*}-1 \\ 
O\left( t^{\nu (2^{*}-1-p)}\right) \quad  & \mathrm{if~}p>2^{*}-1
\end{array}
\right. .  \label{STIME}
\end{equation}
\end{theorem}

Observe that all the cases of (\ref{STIME}) are possible for $2_{\alpha
}<p<2^{*}$ if $N<6$, while only the third case occurs if $N\geq 7$.\medskip

\proof%
By Lemmas \ref{LEM(ptofisso)} and \ref{LEM(pointwise)}, for every $t>0$ we
have 
\[
v\left( t\right) =\frac{B}{t^{\nu }}\left\{ I_{\nu }\left( t\right)
\int_{t}^{+\infty }K\left( s\right) v\left( s\right) ^{p-1}ds+K_{\nu }\left(
t\right) \int_{0}^{t}I\left( s\right) v\left( s\right) ^{p-1}ds\right\} 
\]
and 
\[
v\left( t\right) \leq C_{N,A,\alpha }\left\| v^{\prime }\right\| _{2,\alpha }%
\frac{1}{t^{\nu }}. 
\]
Then, for every $t>0$, one has 
\begin{eqnarray}
v\left( t\right)
&\leq &
C_{N,A,\alpha }^{p-1}\left\| v^{\prime }\right\|
_{2,\alpha }^{p-1}\frac{B}{t^{\nu }}\left\{ I_{\nu }\left( t\right)
\int_{t}^{+\infty }\frac{K\left( s\right) }{s^{\nu \left( p-1\right) }}%
ds+K_{\nu }\left( t\right) \int_{0}^{t}\frac{I\left( s\right) }{s^{\nu
\left( p-1\right) }}ds\right\}   \label{v<w} \\
&=:&
B\,C_{N,A,\alpha }^{p-1}\left\| v^{\prime }\right\| _{2,\alpha}^{p-1}w\left( t\right)   \nonumber
\end{eqnarray}
with obvious definition of $w\left( t\right) $. Note that $w\left( t\right)
\in \mathbb{R}$, by the same reasons used in Remark \ref{RMK(T definito)}.
We now study the behaviour of $w\left( t\right) $ as $t\rightarrow 0^{+}$.

By estimates (\ref{IKzero}) and De L'H\^{o}pital's rule, one obtains 
\[
\int_{0}^{t}\frac{I\left( s\right) }{s^{\nu \left( p-1\right) }}ds\asymp
t^{\nu \left( 2^{*}+1-p\right) } 
\]
and hence, since $t^{-\nu }K_{\nu }\left( t\right) \asymp t^{-2\nu }$, we
have 
\[
\frac{K_{\nu }\left( t\right) }{t^{\nu }}\int_{0}^{t}\frac{I\left( s\right) 
}{s^{\nu \left( p-1\right) }}ds\asymp t^{\nu \left( 2^{*}-1-p\right) }. 
\]
In particular, since both sides are positive, there exists $C_{1}>0$ such
that 
\begin{equation}
\frac{K_{\nu }\left( t\right) }{t^{\nu }}\int_{0}^{t}\frac{I\left( s\right) 
}{s^{\nu \left( p-1\right) }}ds=C_{1}t^{\nu \left( 2^{*}-1-p\right)
}+o\left( t^{\nu \left( 2^{*}-1-p\right) }\right)  \label{w2}
\end{equation}
(one can also check that $C_{1}=1/\left( 2\nu ^{2}(2^{*}+1-p)\right) $).

Assume $2^{*}-1-p>0$. Then 
\[
w\left( t\right) =\frac{I_{\nu }\left( t\right) }{t^{\nu }}\int_{t}^{+\infty
}\frac{K\left( s\right) }{s^{\nu \left( p-1\right) }}ds+\frac{K_{\nu }\left(
t\right) }{t^{\nu }}\int_{0}^{t}I\left( s\right) \frac{I\left( s\right) }{%
s^{\nu \left( p-1\right) }}ds=\frac{I_{\nu }\left( t\right) }{t^{\nu }}%
\int_{t}^{+\infty }\frac{K\left( s\right) }{s^{\nu \left( p-1\right) }}%
ds+o\left( 1\right) . 
\]
Since $K\left( s\right) s^{-\nu \left( p-1\right) }\asymp s^{\frac{N+\alpha 
}{2-\alpha }-\nu p}$ (see (\ref{IKzero})) and 
\[
\frac{N+\alpha }{2-\alpha }-\nu p>\frac{N+\alpha }{2-\alpha }-\nu \left(
2^{*}-1\right) =-1, 
\]
we have $\int_{t}^{+\infty }K\left( s\right) s^{-\nu \left( p-1\right)
}ds\asymp 1$, whence 
\[
\frac{I_{\nu }\left( t\right) }{t^{\nu }}\int_{t}^{+\infty }\frac{K\left(
s\right) }{s^{\nu \left( p-1\right) }}ds\asymp 1 
\]
because $t^{-\nu }I_{\nu }\left( t\right) \asymp 1$. Therefore $w\left(
t\right) \asymp 1$ and the first estimate of (\ref{STIME}) then follows from
(\ref{v<w}).

Now we assume $2^{*}-1-p<0$. Then we have $K\left( s\right) s^{-\nu \left(
p-1\right) }\asymp s^{\frac{N+\alpha }{2-\alpha }-\nu p}$ (see (\ref{IKzero}%
)) and 
\begin{equation}
\frac{N+\alpha }{2-\alpha }-\nu p<\frac{N+\alpha }{2-\alpha }-\nu \left(
2^{*}-1\right) =-1,  \label{diverge}
\end{equation}
so that $\int_{0}^{+\infty }K\left( s\right) s^{-\nu \left( p-1\right)
}ds=+\infty $. By estimates (\ref{IKzero}) and De L'H\^{o}pital's rule, we
get 
\[
\int_{t}^{+\infty }\frac{K\left( s\right) }{s^{\nu \left( p-1\right) }}%
ds\asymp t^{\nu \left( 2^{*}-1-p\right) }, 
\]
which gives 
\[
\frac{I_{\nu }\left( t\right) }{t^{\nu }}\int_{t}^{+\infty }\frac{K\left(
s\right) }{s^{\nu \left( p-1\right) }}ds\asymp t^{\nu \left(
2^{*}-1-p\right) }, 
\]
since $t^{-\nu }I_{\nu }\left( t\right) \asymp 1$. In particular, since both
sides are positive, there exists $C_{2}>0$ such that 
\begin{equation}
\frac{I_{\nu }\left( t\right) }{t^{\nu }}\int_{t}^{+\infty }\frac{K\left(
s\right) }{s^{\nu \left( p-1\right) }}ds=C_{2}t^{\nu \left( 2^{*}-1-p\right)
}+o\left( t^{\nu \left( 2^{*}-1-p\right) }\right)  \label{w1}
\end{equation}
(one can also check that $C_{2}=-1/\left( 2\nu ^{2}(2^{*}-1-p)\right) $).
Therefore, by (\ref{w1}) and (\ref{w2}), we obtain 
\begin{eqnarray*}
w\left( t\right)  &=&\frac{I_{\nu }\left( t\right) }{t^{\nu }}%
\int_{t}^{+\infty }\frac{K\left( s\right) }{s^{\nu \left( p-1\right) }}ds+%
\frac{K_{\nu }\left( t\right) }{t^{\nu }}\int_{0}^{t}I\left( s\right) \frac{%
I\left( s\right) }{s^{\nu \left( p-1\right) }}ds \\
&=&\left( C_{1}+C_{2}\right) t^{\nu \left( 2^{*}-1-p\right) }+o\left( t^{\nu
\left( 2^{*}-1-p\right) }\right) 
\end{eqnarray*}
and thus the third estimate of (\ref{STIME}) follows from (\ref{v<w}).

Finally, we assume $2^{*}-1-p=0$. We get $\int_{0}^{+\infty }K\left(
s\right) s^{-\nu \left( p-1\right) }ds=+\infty $ again (the inequality (\ref
{diverge}) becomes an equality), but the estimates (\ref{IKzero}) and De
L'H\^{o}pital's rule now give 
\[
\int_{t}^{+\infty }\frac{K\left( s\right) }{s^{\nu \left( p-1\right) }}%
ds\asymp \ln t 
\]
and hence 
\[
\frac{I_{\nu }\left( t\right) }{t^{\nu }}\int_{t}^{+\infty }\frac{K\left(
s\right) }{s^{\nu \left( p-1\right) }}ds\asymp \ln t. 
\]
In particular, since the left hand side is positive, there exists $C_{3}>0$
such that 
\[
\frac{I_{\nu }\left( t\right) }{t^{\nu }}\int_{t}^{+\infty }\frac{K\left(
s\right) }{s^{\nu \left( p-1\right) }}ds=-C_{3}\ln t+o\left( \ln t\right) 
\]
(one can also check that $C_{3}=1/(2\nu )$). So, by (\ref{w2}), we conclude
that 
\[
w\left( t\right) =-C_{3}\ln t+o\left( \ln t\right) +C_{1}+o\left( 1\right)
=-C_{3}\ln t+o\left( \ln t\right) 
\]
and therefore the second estimate of (\ref{STIME}) follows from (\ref{v<w}).%
\endproof%
\bigskip

\noindent \textbf{Proof of Theorem \ref{THM(stime)}.}\quad It readily
follows from Theorem \ref{THM(stime_v)}, by the change of variables (\ref
{cambio1})-(\ref{cambio3}).%
\endproof%

\section{Nonexistence result\label{SEC:nonex}}

In this section we assume 
\begin{equation}
0<\alpha <2\quad \textrm{and}\quad \frac{2N}{N-\alpha }<p\leq 2\frac{%
2N-2+\alpha }{2N-2-\alpha }  \label{hp_alfa&p}
\end{equation}
and consider the problem (\ref{Pb_phi}) of the radial solutions $u\left(
x\right) =\phi \left( \left| x\right| \right) $ of (\ref{P}) which belongs
to $L^{p}(\mathbb{R}^{N})$, that is, 
\begin{equation}
\left\{ 
\begin{array}{ll}
\smallskip -\phi ^{\prime \prime }-\dfrac{N-1}{r}\phi ^{\prime }+\dfrac{A}{%
r^{\alpha }}\phi =\phi ^{p-1} & \mathrm{in~}\mathbb{R}_{+} \\ 
\medskip \phi >0 & \mathrm{in~}\mathbb{R}_{+} \\ 
\medskip r^{-\frac{\alpha }{2}}\phi ,\phi ^{\prime }\in L^{2}(\mathbb{R}%
_{+},r^{N-1}dr) &  \\ 
\phi \in L^{p}(\mathbb{R}_{+},r^{N-1}dr) & 
\end{array}
\right. .  \label{nonex_pb}
\end{equation}

We set 
\begin{equation}
\beta :=\frac{\alpha p}{p-2}.  \label{beta:=}
\end{equation}
Notice that, since $\alpha >0$, the second condition of (\ref{hp_alfa&p}) is
equivalent to 
\[
\frac{2N-2+\alpha }{2}\leq \beta <N. 
\]
Moreover, one has 
\[
\beta -2\geq \frac{2N-2+\alpha }{2}-2=\frac{2N-6+\alpha }{2}>0. 
\]

\begin{lemma}
\label{LEM(liminf)}Assume that $\phi $ is a solution of problem (\ref
{nonex_pb}) (with conditions (\ref{hp_alfa&p})). Then 
\[
\lim_{r\rightarrow 0^{+}}r^{\beta -2}\phi \left( r\right) ^{2}=0.
\]
\end{lemma}

\proof%
It follows from Theorem \ref{THM(stime)}. Indeed, if $p>2^{*}-1$, then one
has 
\[
r^{\beta -2}\phi \left( r\right) ^{2}=r^{\beta -2}O\left(
r^{(2^{*}-1-p)\left( N-2\right) }\right) =O\left( r^{(2^{*}-1-p)\left(
N-2\right) +\beta -2}\right) \quad \mathrm{as~}r\rightarrow 0^{+}, 
\]
where 
\begin{eqnarray*}
\left( 2^{*}-1-p\right) \left( N-2\right) +\beta -2 &\geq &\left( 2^{*}-1-2%
\frac{2N-2+\alpha }{2N-2-\alpha }\right) \left( N-2\right) +\frac{%
2N-2+\alpha }{2}-2 \\
&=&\frac{\left( 2-\alpha \right) \left( 6N-6+\alpha \right) }{2\left(
2N-2-\alpha \right) }>0.
\end{eqnarray*}
The other cases are obvious, since $\beta >2$.%
\endproof%
\bigskip

\noindent \textbf{Proof of Theorem \ref{THM(nonex)}.}\quad For the sake of
contradiction, we assume the $\phi $ is a solution of problem (\ref{nonex_pb}%
) (with conditions (\ref{hp_alfa&p})). Rewriting the equation of (\ref
{nonex_pb}) in the following form 
\[
r^{1-N}\left( r^{N-1}\phi ^{\prime }\right) ^{\prime }-\dfrac{A}{r^{\alpha }}%
\phi +\phi ^{p-1}=0\quad \mathrm{in~}\mathbb{R}_{+} 
\]
and testing it with $r^{\beta -1}\phi \left( r\right) $ on an arbitrary
interval $\left[ a,b\right] \subset \mathbb{R}_{+}$, we get 
\begin{equation}
\int_{a}^{b}\left( \left( r^{N-1}\phi ^{\prime }\right) ^{\prime }r^{\beta
-N}\phi -Ar^{\beta -1-\alpha }\phi ^{2}+r^{\beta -1}\phi ^{p}\right) dr=0.
\label{eq_test}
\end{equation}
Integrating by parts twice, one finds that 
\begin{eqnarray*}
\int_{a}^{b}\left( r^{N-1}\phi ^{\prime }\right) ^{\prime }r^{\beta -N}\phi
\,dr &=&\left[ r^{\beta -1}\phi ^{\prime }\phi \right]
_{a}^{b}-\int_{a}^{b}\left( r^{\beta -1}\left( \phi ^{\prime }\right)
^{2}+\left( \beta -N\right) r^{\beta -2}\phi \phi ^{\prime }\right) dr \\
&=&\left[ r^{\beta -1}\phi ^{\prime }\phi \right]
_{a}^{b}-\int_{a}^{b}r^{\beta -1}\left( \phi ^{\prime }\right) ^{2}dr+ \\
&&-\frac{\beta -N}{2}\left[ r^{\beta -2}\phi ^{2}\right] _{a}^{b}+\frac{%
\left( \beta -N\right) \left( \beta -2\right) }{2}\int_{a}^{b}r^{\beta
-3}\phi ^{2}dr,
\end{eqnarray*}
so that, plugging into (\ref{eq_test}), we obtain 
\[
\left[ r^{\beta -1}\phi ^{\prime }\phi -\frac{\beta -N}{2}r^{\beta -2}\phi
^{2}\right] _{a}^{b}-\int_{a}^{b}r^{\beta -1}\left( \phi ^{\prime }\right)
^{2}dr+\frac{\left( \beta -N\right) \left( \beta -2\right) }{2}%
\int_{a}^{b}r^{\beta -3}\phi ^{2}dr+ 
\]
\begin{equation}
-A\int_{a}^{b}r^{\beta -1-\alpha }\phi ^{2}dr+\int_{a}^{b}r^{\beta -1}\phi
^{p}dr=0.  \label{i}
\end{equation}

We now define 
\begin{equation}
E\left( r\right) :=\frac{1}{2}\phi ^{\prime }\left( r\right) ^{2}-\frac{1}{2}%
\dfrac{A}{r^{\alpha }}\phi \left( r\right) ^{2}+\frac{1}{p}\phi \left(
r\right) ^{p}\quad \mathrm{and}\quad E_{\beta }\left( r\right) :=r^{\beta
}E\left( r\right) \quad \mathrm{for~all~}r>0.  \label{E:=}
\end{equation}
Taking the derivative of $E$ and using the equation, we get 
\begin{eqnarray*}
E^{\prime }\left( r\right) &=&\phi ^{\prime \prime }\phi ^{\prime }+\frac{%
\alpha }{2}\dfrac{A}{r^{\alpha +1}}\phi ^{2}-\dfrac{A}{r^{\alpha }}\phi \phi
^{\prime }+\phi ^{p-1}\phi ^{\prime } \\
&=&\left( -\dfrac{N-1}{r}\phi ^{\prime }+\dfrac{A}{r^{\alpha }}\phi -\phi
^{p-1}\right) \phi ^{\prime }+\frac{\alpha }{2}\dfrac{A}{r^{\alpha +1}}\phi
^{2}-\dfrac{A}{r^{\alpha }}\phi \phi ^{\prime }+\phi ^{p-1}\phi ^{\prime } \\
&=&-\dfrac{N-1}{r}\left( \phi ^{\prime }\right) ^{2}+\frac{\alpha }{2}\dfrac{%
A}{r^{\alpha +1}}\phi ^{2}
\end{eqnarray*}
and hence 
\begin{eqnarray}
E_{\beta }\left( b\right) -E_{\beta }\left( a\right) &=&\int_{a}^{b}\left(
\beta r^{\beta -1}E\left( r\right) +r^{\beta }E^{\prime }\left( r\right)
\right) dr  \nonumber \\
&=&\left( \frac{\beta }{2}-N+1\right) \int_{a}^{b}r^{\beta -1}\left( \phi
^{\prime }\right) ^{2}dr+
\label{ii} \\
&& +\frac{A\left( \alpha -\beta \right) }{2}\int_{a}^{b}r^{\beta -\alpha -1}\phi ^{2}dr
+\frac{\beta }{p}\int_{a}^{b}r^{\beta -1}\phi ^{p}dr.  \nonumber
\end{eqnarray}

Multiplying (\ref{i}) by $\beta /p$ and adding side by side to (\ref{ii}),
we finally obtain 
\[
\left( \frac{\beta }{p}+\frac{\beta }{2}-N+1\right) \int_{a}^{b}r^{\beta
-1}\left( \phi ^{\prime }\right) ^{2}dr+A\left( \frac{\alpha -\beta }{2}+%
\frac{\beta }{p}\right) \int_{a}^{b}r^{\beta -\alpha -1}\phi ^{2}dr 
\]
\begin{equation}
+\frac{\beta \left( N-\beta \right) \left( \beta -2\right) }{2p}%
\int_{a}^{b}r^{\beta -3}\phi ^{2}dr=\frac{\beta }{p}\left[ r^{\beta -1}\phi
^{\prime }\phi -\frac{\beta -N}{2}r^{\beta -2}\phi ^{2}\right]
_{a}^{b}+E_{\beta }\left( b\right) -E_{\beta }\left( a\right) ,  \label{i+ii}
\end{equation}
where the second term of the left hand side actually vanishes, since $%
(\alpha -\beta )/2+\beta /p=0$ thanks to the definition (\ref{beta:=}) of $%
\beta $.

We now use the integrability properties (\ref{nonex_pb}) of $\phi $ and $%
\phi ^{\prime }$. Since $\beta <N$, we have 
\[
\beta -1-N+\frac{\alpha }{2}<\frac{\alpha }{2}-1<0\quad \textrm{and}\quad
\beta -3<N-3<N-1-\alpha , 
\]
so that 
\begin{eqnarray*}
\int_{1}^{+\infty }r^{\beta -2}\left| \phi ^{\prime }\right| \phi \,dr
&=&\int_{1}^{+\infty }r^{\frac{N-1}{2}}\left| \phi ^{\prime }\right| \,r^{%
\frac{N-1-\alpha }{2}}\phi \,r^{\beta -1-N+\frac{\alpha }{2}}dr\\
&\leq&
\int_{1}^{+\infty }r^{\frac{N-1}{2}}\left| \phi ^{\prime }\right| \,r^{\frac{N-1-\alpha }{2}}\phi \,dr \\
&\leq &
\left( \int_{1}^{+\infty }r^{N-1}\left( \phi ^{\prime }\right)
^{2}dr\right) ^{1/2}\left( \int_{1}^{+\infty }r^{N-1-\alpha }\phi
^{2}dr\right) ^{1/2}<\infty
\end{eqnarray*}
and 
\begin{eqnarray*}
&&\int_{1}^{+\infty }\left( r^{\beta -3}\phi ^{2}+r^{\beta -1}\left( \phi
^{\prime }\right) ^{2}+r^{\beta -1-\alpha }\phi ^{2}+r^{\beta -1}\phi
^{p}\right) dr \\
&\leq &\int_{1}^{+\infty }\left( r^{\beta -1-\alpha }\phi ^{2}+r^{N-1}\left(
\phi ^{\prime }\right) ^{2}+r^{N-1-\alpha }\phi ^{2}+r^{N-1}\phi ^{p}\right)
dr<\infty .
\end{eqnarray*}
This implies 
\[
\liminf_{r\rightarrow +\infty }\left( r^{\beta -1}\left| \phi ^{\prime
}\right| \phi +r^{\beta -2}\phi ^{2}+r^{\beta }\left( \phi ^{\prime }\right)
^{2}+r^{\beta -\alpha }\phi ^{2}+r^{\beta }\phi ^{p}\right) =0 
\]
and thus there exists a sequence $b_{n}\rightarrow +\infty $ such that 
\[
\lim_{n\rightarrow \infty }b_{n}^{\beta -1}\phi ^{\prime }\left(
b_{n}\right) \phi \left( b_{n}\right) =\lim_{n\rightarrow \infty
}b_{n}^{\beta -2}\phi \left( b_{n}\right) ^{2}=\lim_{n\rightarrow \infty
}E_{\beta }\left( b_{n}\right) =0. 
\]
Evaluating (\ref{i+ii}) with $b=b_{n}$ and passing to the limit, we find 
\begin{eqnarray}
\gamma _{1}\int_{a}^{+\infty }r^{\beta -1}\left( \phi ^{\prime }\right)
^{2}dr+\gamma _{2}\int_{a}^{+\infty }r^{\beta -3}\phi ^{2}dr &=&-\frac{\beta 
}{p}a^{\beta -1}\phi ^{\prime }\left( a\right) \phi \left( a\right) +   \label{IDENTITA'}\\
&&-\frac{\beta \left( N-\beta \right) }{2p}a^{\beta -2}\phi \left( a\right)
^{2}-E_{\beta }\left( a\right) ,  \nonumber
\end{eqnarray}
where 
\[
\gamma _{1}:=\frac{\beta }{p}+\frac{\beta }{2}-N+1,\quad \gamma _{2}:=\frac{%
\beta \left( N-\beta \right) \left( \beta -2\right) }{2p}. 
\]
We study the two sides of identity (\ref{IDENTITA'}) separately. We have 
\[
\gamma _{1}=\frac{2N-2-\alpha }{2\left( p-2\right) }\left( 2\frac{%
2N-2+\alpha }{2N-2-\alpha }-p\right) \geq 0 
\]
and $\gamma _{2}>0$ (recall that $\beta <N$ and $\beta >2$). As a
consequence, since $\phi >0$, there exist two constants $a_{0},\gamma _{0}>0$
such that 
\begin{equation}
\forall a\leq a_{0},\quad \gamma _{1}\int_{a}^{+\infty }r^{\beta -1}\left(
\phi ^{\prime }\right) ^{2}dr+\gamma _{2}\int_{a}^{+\infty }r^{\beta -3}\phi
^{2}dr\geq \gamma _{0}>0.  \label{1°m}
\end{equation}
On the other hand, Lemma \ref{LEM(liminf)} assures that 
\[
\lim_{a\rightarrow 0^{+}}a^{\beta -2}\phi \left( a\right) ^{2}=0, 
\]
which also gives 
\[
\lim_{a\rightarrow 0^{+}}a^{\beta -\alpha }\phi \left( a\right)
^{2}=\lim_{a\rightarrow 0^{+}}a^{2-\alpha }a^{\beta -2}\phi \left( a\right)
^{2}=0. 
\]
Therefore, briefly denoting the right hand side of (\ref{IDENTITA'}) by $%
F\left( a\right) $ and substituting the definitions (\ref{E:=}), we infer
that 
\begin{eqnarray}
F\left( a\right) 
&=&
-\frac{\beta }{p}a^{\beta -1}\phi ^{\prime }\left(
a\right) \phi \left( a\right) -\frac{\beta \left( N-\beta \right) }{2p}%
a^{\beta -2}\phi \left( a\right) ^{2}-\frac{1}{2}a^{\beta }\phi ^{\prime
}\left( a\right) ^{2}+\frac{A}{2}a^{\beta -\alpha }\phi \left( a\right) ^{2}+  \nonumber \\
&&
-\frac{1}{p}a^{\beta }\phi \left( a\right) ^{p}
-\frac{1}{2}\frac{\beta ^{2}}{p^{2}}a^{\beta -2}\phi \left( a\right) ^{2}+%
\frac{1}{2}\frac{\beta ^{2}}{p^{2}}a^{\beta -2}\phi \left( a\right) ^{2} 
\nonumber \\
&=&-\frac{1}{2}a^{\beta -2}\left( a\phi ^{\prime }\left( a\right) +\frac{%
\beta }{p}\phi \left( a\right) \right) ^{2}-\frac{1}{p}a^{\beta }\phi \left(
a\right) ^{p}+o\left( 1\right) _{a\rightarrow 0^{+}}\leq o\left( 1\right)
_{a\rightarrow 0^{+}}.  \label{2°m}
\end{eqnarray}
So, from (\ref{IDENTITA'}), (\ref{1°m}) and (\ref{2°m}) it follows that $%
\forall a\leq a_{0}$ one has $0<\gamma _{0}\leq F\left( a\right) \leq
o\left( 1\right) _{a\rightarrow 0^{+}}$, which is a contradiction.%
\endproof%

\section{Appendix}

This Appendix is devoted to a summary of the most useful properties of the
Bessel functions used in the paper. For a complete treatment, we refer the
reader to \cite{Niki-Uva}, \cite{Temme} and \cite{Watson}.

For every $\nu \in \mathbb{R}$ and $t\in \mathbb{R}_{+}$, the \emph{modified
Bessel function of the first kind of order} $\nu $ is defined as 
\[
I_{\nu }\left( t\right) =\left( \frac{t}{2}\right) ^{\nu }\sum_{k=0}^{\infty
}\frac{1}{k!}\frac{1}{\Gamma \left( \nu +k+1\right) }\left( \frac{t}{2}%
\right) ^{2k}, 
\]
where $\Gamma $ is the usual Gamma function and $1/\Gamma \left( -n\right)
=0 $ for $n\in \mathbb{N}$. The \emph{modified Bessel function of the second
kind of order} $\nu $ (also known as \emph{Macdonald's function}) is defined
as 
\[
K_{\nu }\left( t\right) =\frac{\pi }{2}\frac{I_{-\nu }\left( t\right)
-I_{\nu }\left( t\right) }{\sin \left( \pi \nu \right) }\quad \mathrm{if~}%
\nu \notin \mathbb{Z} 
\]
and $K_{n}\left( t\right) =\lim_{\nu \rightarrow n}K_{\nu }\left( t\right) $
if $n\in \mathbb{Z}$. These functions are linearly independent real solutions
of the \emph{modified Bessel equation in }$\mathbb{R}_{+}$, namely, 
\[
-u^{\prime \prime }-\dfrac{1}{t}u^{\prime }+\left( 1+\dfrac{\nu ^{2}}{t^{2}}%
\right) u=0\quad \mathrm{in\ }\mathbb{R}_{+}, 
\]
and satisfy the following identities on $\mathbb{R}_{+}$:

\begin{itemize}
\item  $\displaystyle K_{\nu +1}\left( t\right) I_{\nu }\left( t\right)
+I_{\nu +1}\left( t\right) K_{\nu }\left( t\right) =\frac{1}{t};$

\item  $\displaystyle I_{\nu }^{\prime }\left( t\right) -\frac{\nu }{t}%
I_{\nu }\left( t\right) =I_{\nu +1}\left( t\right) ;$

\item  $\displaystyle K_{\nu }^{\prime }\left( t\right) -\frac{\nu }{t}%
K_{\nu }\left( t\right) =-K_{\nu +1}\left( t\right) .$
\end{itemize}

For every $\nu >0$, both $I_{\nu }$ and $K_{\nu }$ are strictly positive on $%
\mathbb{R}_{+}$ and the following asymptotic estimates hold: 
\[
\begin{array}{llll}
I_{\nu }\left( t\right) \sim \dfrac{1}{2^{\nu }\Gamma \left( \nu +1\right) }%
t^{\nu }\quad & \mathrm{and}\quad & K_{\nu }\left( t\right) \sim \dfrac{%
\Gamma \left( \nu \right) }{2^{1-\nu }}t^{-\nu } & \mathrm{as~}t\rightarrow
0^{+},\bigskip \\ 
I_{\nu }\left( t\right) \sim \sqrt{\dfrac{1}{2\pi }}\dfrac{e^{t}}{\sqrt{t}}
& \mathrm{and} & K_{\nu }\left( t\right) \sim \sqrt{\dfrac{\pi }{2}}\dfrac{%
e^{-t}}{\sqrt{t}}\quad & \mathrm{as~}t\rightarrow +\infty .
\end{array}
\]


\begin{thebibliography}{99}
\bibitem{BBR 1}  M. Badiale, V. Benci, S. Rolando, Solitary waves: physical
aspects and mathematical results, Rend. Sem. Math. Univ. Pol. Torino 62
(2004), 107-154.

\bibitem{BGR}  M. Badiale, M. Guida, S. Rolando, Elliptic equations with
decaying cylindrical potentials and power-type nonlinearities, Adv.
Differential Equations 12 (2007), 1321-1362.

\bibitem{BGRrmks}  M. Badiale, M. Guida, S. Rolando, Some general existence
results for nonlinear elliptic equations with potentials, work in progress.

\bibitem{BPR}  M. Badiale, L. Pisani, S. Rolando, Sum of weighted Lebesgue
spaces and nonlinear elliptic equations, NoDEA, Nonlinear Differ. Equ. Appl.
18 (2011), 369-405.

\bibitem{BRpow}  M. Badiale, S. Rolando, A note on nonlinear elliptic
problems with singular potentials, Rend. Lincei Mat. Appl. 16 (2006), 1-13.

\bibitem{BRnorm}  M. Badiale, S. Rolando, Vortices with prescribed $L^{2}$\
norm in the nonlinear wave equation, Adv. Nonlinear Stud. \textbf{8} (2008),
817-842.

\bibitem{BR TMA}  M. Badiale, S. Rolando, Nonlinear elliptic equations with
subhomogeneous potentials, Nonlinear Anal. \textbf{72} (2010), 602-617.

\bibitem{BRcharge}  M. Badiale, S. Rolando, A note on vortices with
prescribed charge, Adv. Nonlinear Stud. \textbf{12} (2012), 703-716.

\bibitem{BF sw intro}  V. Benci, D. Fortunato, Solitary waves in the
nonlinear wave equation and in gauge theories, J. Fixed Point Theory Appl. 1
(2007), 61-86.

\bibitem{Beres-Lions}  H. Berestycki, P.L. Lions, Nonlinear Scalar Field
Equations. I - Existence of a Ground State, Arch. Rational Mech. Anal. 82
(1983), 313-345.

\bibitem{Co-Cr-Par}  M. Conti, S. Crotti, D. Pardo, On the existence of
positive solutions for a class of singular elliptic equations, Adv.
Differential Equations 3 (1998), 111-132.

\bibitem{Co-Terr-Verz}  M. Conti, S. Terracini, G. Verzini, Nodal solutions
to a class of nonstandard superlinear equations on $\mathbb{R}^{N}$, Adv.
Differential Equations 7 (2002), 297-318.

\bibitem{G-R}  M. Guida, S. Rolando, On the asymptotic behaviour of weak
solutions to nonlinear elliptic equations with potential, work in progress.

\bibitem{Lieb-Loss}  E.H. Lieb, M. Loss, Analysis. Second edition, American
Mathematical Soc., 2001.

\bibitem{Niki-Uva}  F. Nikiforov, V. Uvarov, Special functions of
Mathematical Physics. A unified introduction with applications, Birkh\"{a}%
user, 1988.

\bibitem{Pucci-Serrin}  P. Pucci, J. Serrin, A general variational identity,
Indiana Univ. Math. J. 35 (1986), 681-703.

\bibitem{THESIS}  S. Rolando, Nonlinear elliptic equations with singular
symmetric potentials, PhD Thesis, Dipartimento di Matematica, Universit\`{a}
degli Studi di Torino, 2006.\newline
(www2.dm.unito.it/paginepersonali/rolando)

\bibitem{Strauss L.Notes}  W.A. Strauss, Nonlinear invariant wave equations,
Lecture Notes in Physics, vol. 23, Springer, 1978.

\bibitem{Su-Wang-Will 2}  J. Su, Z.Q. Wang, M. Willem, Nonlinear Schr\"{o}%
dinger equations with unbounded and decaying potentials, Commun. Contemp.
Math. 9 (2007), 571-583.

\bibitem{Su-Wang-Will p}  J. Su, Z.Q. Wang, M. Willem, Weighted Sobolev
embedding with unbounded and decaying radial potentials, J. Differential
Equations 238 (2007), 201-219.

\bibitem{Temme}  N. Temme, Special functions. An introduction to the
classical functions of Mathematical Physics, John Wiley \& Sons, 1996.

\bibitem{Terracini}  S. Terracini, On positive entire solutions to a class
of equations with singular coefficient and critical exponent, Adv.
Differential Equations 1 (1996), 241-264.

\bibitem{Watson}  G.N. Watson, A treatise on the theory of Bessel functions,
Cambridge University Press, 1952.

\bibitem{YangY}  Y. Yang, Solitons in field theory and nonlinear analysis,
Springer, 2001.\bigskip \bigskip 
\end{thebibliography}
\end{document}